\documentclass[11pt]{article}

\usepackage[a4paper,margin=1in]{geometry}
\usepackage{amsmath,amssymb,amsfonts}
\usepackage{braket}
\usepackage{graphicx}
\usepackage{subcaption}
\usepackage{xcolor} 
\usepackage{hyperref}
\usepackage{cite}
\usepackage{tikz}

\title{From quantum to quantum-inspired: the `LogQ' algorithm as a non-linear continuous relaxation of variables method}

\author{
Jérémie Messud\\
\textit{TotalEnergies}\\
2, place Jean Millier\\
92078 Paris La Défense Cedex, France\\
\texttt{jeremie.messud@totalenergies.com}
\and
Yagnik Chatterjee\\
\textit{TotalEnergies}\\
2, place Jean Millier\\
92078 Paris La Défense Cedex, France\\
\texttt{yagnik.chatterjee@totalenergies.com}
}

\begin{document}

\maketitle

\begin{abstract}
The LogQ algorithm encodes Quadratic Unconstrained Binary Optimization (QUBO) problems,
which are often encountered in the industry (portfolio optimization, fleet optimization, charging stations, etc.).
It was developed within the framework of quantum computing, designed as a pragmatic approach to quantum combinatorial optimization that drastically reduces the number of required qubits and quantum circuit depth.
While LogQ has recently been made compliant with gradient-inspired methods, greatly improving parameter optimization efficiency, it still faced hurdles regarding Pauli decomposition and measurement overhead.
We here demonstrate that LogQ can be fully reformulated within a classical framework, which effectively eliminates the need for Pauli decomposition and  bypasses the measurement challenges altogether.
This finally leads to a classical heuristic based on a non-linear continuous relaxation of variables and is, to the best of our knowledge, novel.
The LogQ story illustrates how quantum computing can inspire classical algorithms, leading to so-called "quantum-inspired" methods.
\end{abstract}


\section{Introduction}

Quadratic Unconstrained Binary Optimization (QUBO) problems have been widely studied \cite{qubotutorial, date2021qubo, calude2017qubo, papalitsas2019qubo} and appear in many industrial applications (portfolio optimization, fleet optimization, charging stations, etc.).
Their worst-case complexity is NP-hard \cite{karp, surveypoly}, which has motivated the development of heuristic methods capable of delivering good approximations in reasonable time.
The Branch \& Bound (B\&B) algorithm, which among others relies on linear relaxation of variables, usually represents the standard baseline\cite{land1960automatic}.
However, problems of industrial interest can still be hard for such classical heuristics \cite{koch2025quantumoptimizationbenchmarkinglibrary}.

Quantum computing (QC) based heuristics, especially hybrid quantum-classical algorithms such as QAOA \cite{qaoa,qaoa2}, VQE \cite{peruzzo2014variational}, quantum annealing \cite{Farhi2001QuantumAnnealing} and LogQ \cite{rancic, chatterjee2023solving, chatterjee2023hybrid, neasqc, chatterjee2024quantum}, have also been investigated for such tasks.
In contrast to QAOA, which uses computational basis state encoding, LogQ relies on amplitude-encoding, which leads to an exponential reduction in the number of qubits 
and a quadratic reduction in quantum circuit depth, that can be attractive for industrial-sized instances.

The original LogQ formulation was not compliant with gradient-inspired optimization of its parameters (due to a vanishing gradient issue), and therefore required evolutionary or global optimization methods, potentially hampering its efficiency.
A recent evolution of LogQ addresses this issue and yields a gradient-inspired methods compliant landscape \cite{chatterjee2025logq} (when we mention gradient-inspired methods we essentially think perturbed-gradient methods).
However, a challenge remains: LogQ requires a Pauli decomposition of the QUBO operator, which is itself costly, and many measurements, which can hinder quantum advantage.

In this paper, we overcome this challenge by showing that LogQ can in fact naturally be reformulated within a classical framework without Pauli decomposition and several measurements.
This leads to a clear classical interpretation: we obtain a continuous relaxation of variables method which operates non-linearly on `phases' and is new to our knowledge. 
This provides an illustration of how quantum constructions can inspire classical algorithms, leading to `quantum-inspired' methods.

\section{Spin-QUBO problems and linear relaxation of variables}

We adopt spin variables $s_i \in \{-1,1\}$ rather than binary variables $x_i \in \{0,1\}$ to maintain consistency with the quantum formalism.
They are related by a simple relation:
    $s_i = 2x_i-1$.
%
%
The spin-QUBO (sQUBO) version of the QUBO problem 
is, for $n$ variables \cite{chatterjee2024quantum}:
\begin{align}
\label{eq:sQUBO}
s^*=\underset{s \in \{-1,1\}^{n}}{\text{argmin}}-\frac{1}{2}\sum_{i,j=0}^{n-1} s_i Q_{ij} s_j
.
\end{align}
The $Q_{ij}$ here are the sQUBO matrix coefficients.
%
A potential heuristic to solve this  problem is obtained by linear relaxation, i.e. taking:
\begin{align}
\label{eq:lin-relax}
s_i\rightarrow s^{(lin)}_i\in[-1,1],
\end{align}
with:
\begin{align}
\label{eq:lin-relax2}
\hat{s}_i^*=\begin{cases} 
-1 & \text{if } s^{(lin)*}_i < 0 \\
+1 & \text{if } s^{(lin)*}_i \ge 0 
\end{cases}\quad \leadsto\ s^*_i.
\end{align}
The direct usage of such a method usually leads to poor results. One of the reasons for this is that
$|s^{(lin)}_i|\approx 1$ is not constrained, and thus $|s^{(lin)}_i|\in[0,1]$ can be far from $1$ and the thresholding in (\ref{eq:lin-relax2}) can be very `aggressive'.
A much better option is to embed the linear relaxation within B\&B \cite{land1960automatic},
which usually represents the standard baseline (implemented using the binary version instead of the spin-one but again both are straightforwardly related).
However, problems of industrial interest can still be hard for such classical heuristics \cite{koch2025quantumoptimizationbenchmarkinglibrary}.

\section{Quantum LogQ formulation}
\label{sec:Laplacian}

QC based heuristics, especially hybrid quantum-classical algorithms, have also been investigated for sQUBO tasks.
The LogQ quantum algorithm relies on amplitude-encoding and is formulated as:
\begin{equation}
\label{eqn2}
    \theta^*=\underset{\theta\in \mathbb{R}^n}{\text{argmin }} C(\theta)
    \begin{cases}
          C(\theta)=-2^{N-2}\bra{\Psi(\theta)}\hat{L}\ket{\Psi(\theta)}  \\
          \ket{\Psi(\theta)} = \frac{1}{n}\sum_{i=0}^{n-1}f(\theta_i)\ket{i},
     \end{cases}
\end{equation}
where:
\begin{itemize}
    \item 
$\ket{\Psi(\theta)}$ is a $N$-qubit state with $N=\lceil\log_2n \rceil$ ($\lceil.\rceil$ stands for the ceiling function) - note that logarithmically less qubits than the number of variables are required, hence the LogQ name,
    \item 
$\{\ket{i}, i=0\dots n-1\}$ are the computational basis states,
    \item 
$\theta=(\theta_0\dots\theta_{n-1})$
represent $n$ parameters (equal to the number of spin variables) to optimize classically,
    \item 
The operator $\hat L$ is defined by:
\begin{equation}
\label{L_decompo}
\hat{L}=\dfrac{1}{n}\sum\limits_{k=1}^{n^2}\mathrm{Tr}(J_k Q)J_k,
\end{equation}
where the $J_k$ represent the $n^2=4^N$ possible tensor products of $N$ Pauli matrices and the identity matrix. This Pauli decomposition of the sQUBO operator requires considering at most $(n^2+n)/2$ components (as the matrix $Q$ is symmetric) as well as the same number of measurements on the Pauli basis.
\end{itemize}

The computational basis state amplitudes $f(\theta_i)$ in (\ref{eqn2}) are defined by:
\begin{align}
\label{eq:f}
f(\theta_i) = e^{-i\pi R(\theta_i)},
\end{align}
where the scalar function $R(\cdot)$ is constrained to satisfy:
\begin{align}
\label{con_1}
R(\cdot)\in\mathbb{R},\quad
R(\theta_i^*)=\text{$0$ or $1$ (only) once }\theta\text{ parameters optimized.}
\end{align}
Then, the sQUBO solution can be encoded in $R(\theta_i^*)$:
\begin{align}
\label{sol}
s^*_i=\begin{cases} 
-1 & \text{if } R(\theta_i^*)=0 \\
+1 & \text{if } R(\theta_i^*)=1 
\end{cases}.
\end{align}
Defining an $R(\cdot)$ that satisfies (\ref{con_1}) constitutes the sole degree of freedom in the LogQ scheme and  provides a foundation for developing efficient parameterizations and heuristics.

The conditions in \eqref{con_1} are not constructive, which recently led us to develop more constructive sufficient conditions on $R(\cdot)$ (the first two ones in next equation).
Also, the original proposals in \cite{chatterjee2023solving, rancic} were not gradient-optimization compliant and required evolutionary or even global optimization methods, which are costly.
We developed conditions that allow us to train LogQ (i.e. $\theta$) using gradient-inspired methods \cite{chatterjee2025logq}
 (the last two ones in next equation).
\begin{align}
\label{con_2}
& R(\theta_i)\in[0,1]\ \text{\& differentiable almost everywhere,}\\
& R'(\theta_i)=0 \Rightarrow R(\theta_i)\in\{0,1\},
\nonumber\\
& \text{There exist large intervals such that }|R'(\theta_i)|\ \text{is not small,}
\nonumber\\
& \text{Only `small barrier' local minima exist in $C(\theta)$}.
\nonumber
\end{align}
The final two conditions are intentionally flexible, encouraging the use of perturbed-gradient methods, which can enable us to escape 'small' plateaus (areas where the gradient tends to vanish) and 'small barrier' local minima in the cost-function valley by 'jumping' above these barriers. These flexible constraints lead in practice to heuristic schemes.
The sigmoid combination-based $R(\cdot)$ defined in \cite{chatterjee2025logq}, reasonably satisfying (\ref{con_1}) and (\ref{con_2}), demonstrated robustness.
\begin{figure}[!h]
\begin{center}
  \includegraphics[scale=0.43]{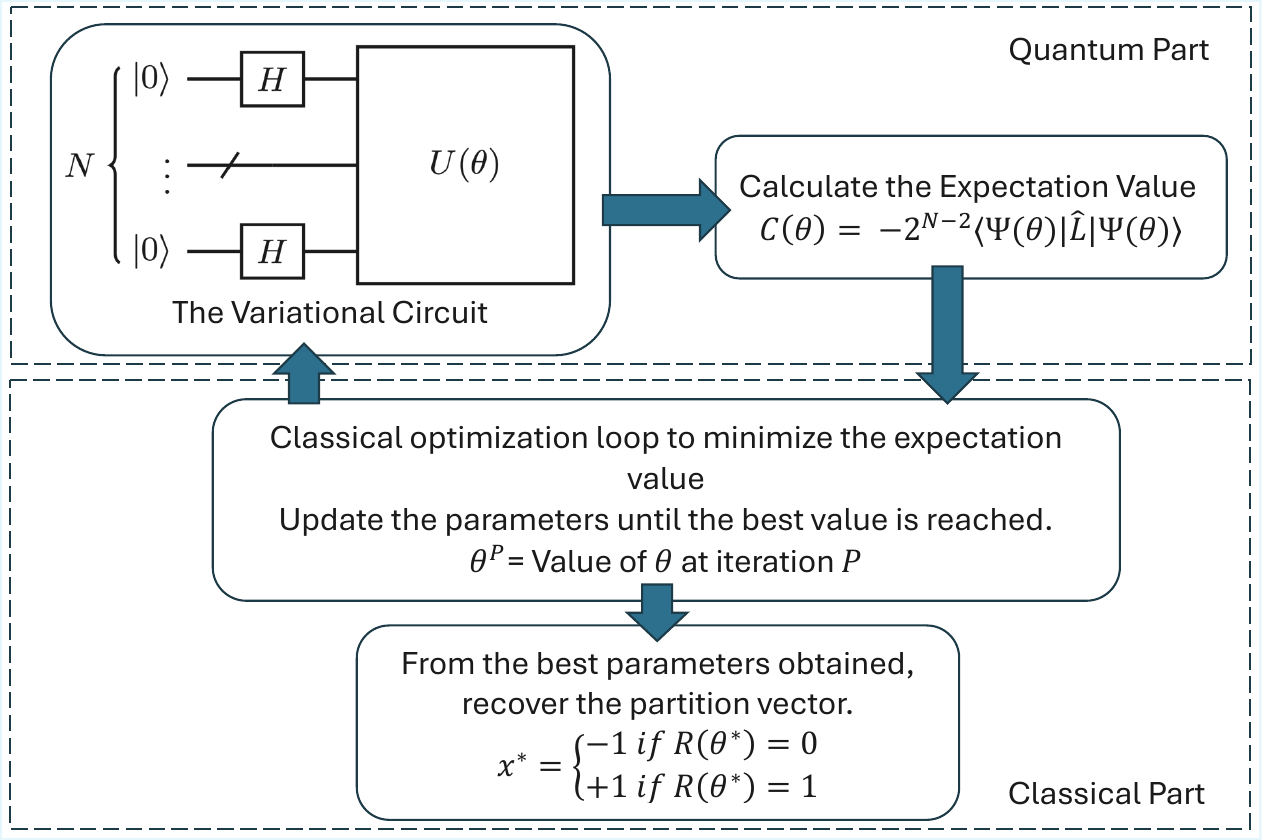}
  \caption{Diagrammatic representation of the hybrid quantum-classical LogQ algorithm.}
  \label{algprocess}
\end{center}
\end{figure}

Finally, note that the $N$-qubit state in (\ref{eqn2}) can equivalently be written as:
\begin{align}
\ket{\Psi(\theta)}= U(\theta)H^{\otimes N}\ket{0}^{\otimes N},
\end{align}
where $U(\theta)$ is a complex and unitary operator (the essence of quantum computation) and is useful to implement LogQ as a quantum circuit on a quantum processing unit (QPU).
LogQ is originally a hybrid quantum-classical algorithm as the optimization
of the $\theta$ parameters is performed using classical methods. The LogQ workflow is shown in Fig.~\ref{algprocess}.

Various works have discussed the potential advantage of quantum LogQ over QAOA on the MaxCut problem \cite{chatterjee2023solving}.
We also did tests on portfolio optimization cases that confirmed the potential advantage of quantum LogQ over QAOA.
However, a challenge remains: LogQ requires a Pauli decomposition of the QUBO operator, eq. (\ref{L_decompo}), which is itself costly and leads to many measurements in the LogQ workflow in Fig. \ref{algprocess}, which can hinders quantum advantage and remains the main challenge of the quantum LogQ algorithm. When trying to overcome this challenge, we realized that LogQ can in fact be fully reformulated  classically and leads to a new continuous relaxation of variables method.

\section{The classical LogQ formulation as a new continuous relaxation of variables method}

The linear relaxation of variables in (\ref{eq:lin-relax}) does not constrain $|s_i|\approx 1$, which hampers its efficiency when used as such as already discussed. 
On the real line, adding this constraint would lead to a problem equivalent to the original one and thus as hard.
The idea is to try to formulate the constraint on the complex plane, which is related to the thought process underlying LogQ as we will discuss. Instead of the form in (\ref{eq:sQUBO})-(\ref{eq:lin-relax}), we take the following form:
\begin{align}
\label{eq:QUBO}
& s^{(LogQ)}_i\leadsto f(\theta_i)\in\mathbb{C},\quad
|f(\theta_i)|=1,\quad
f(\theta_i^*)=\text{$-1$ or $1$ (only) once }\theta\text{ parameters optimized},\nonumber\\
& \theta^*=\underset{\theta \in \mathbb{R}^n}{\text{argmin}}-\frac{1}{2}\sum_{i,j=1}^n f^\dagger(\theta_i) Q_{ij} f(\theta_j)
.
\end{align}
This is no longer a linear relaxation leveraging a real interval, but a continuous relaxation leveraging the complex plane (recall that there are as many continuous variables $\theta_i$ to optimize as the number of original discrete variables).
After optimization, the scheme leads to:
\begin{align}
f(\theta_i^*) \leadsto\ s^*_i,
\end{align}
which does not require any thresholding like in \eqref{eq:lin-relax2}.
How to define $f(\cdot)$? A straightforward choice that satisfies the first line of (\ref{eq:QUBO}) is given by eqs. (\ref{eq:f})-(\ref{con_1}), leading to the decoding rule in (\ref{sol}) after optimization.
The additional conditions in (\ref{con_2}) make full sense in this classical formulation as well (and thus the sigmoid combination-based $R(\cdot)$ studied in \cite{chatterjee2025logq}), for the same reasons than the ones explicated above.
These flexible constraints lead in practice to heuristic schemes.
The final scheme is a classical reformulation of LogQ, equivalent to the quantum LogQ scheme but overcoming the Pauli decomposition and measurement challenge.

Interestingly, these choices make the objective function always real (not only at convergence) \cite{chatterjee2025logq}:
\begin{align}
\label{eq:QUBO2}
\frac{1}{2}\sum_{i,j=1}^n f^\dagger(\theta_i) Q_{ij} f(\theta_j)
=
\sum_{i>j=1}^n \cos\!\big(\pi(R(\theta_i)-R(\theta_j))\big) Q_{ij} + ctnt
.
\end{align}

The non-linear coupling between $R(\theta_i)$ and $R(\theta_j)$ within the cosine term echoes the use of the complex plane under a unit-modulus constraint, $|f(\theta_i)|=1$. This reflects the underlying quantum inspiration, as quantum circuits are inherently unitary and operate over complex Hilbert spaces.
Compared to linear relaxation of variables, the only change here is to replace products by cosine couplings:
\begin{align}
\label{eq:diff}
s^{(lin)}_i s^{(lin)}_j\in[-1,1]
\quad\leftrightarrow\quad
\cos\!\big(\pi(R(\theta_i)-R(\theta_j))\big)\in[-1,1].
\end{align}
LogQ can thus be viewed as a continuous but non-linear relaxation of variables method, operating on phases.
As such, LogQ leads to a classical and apparently new heuristic for combinatorial optimization while benefiting from the previously developed conditions in (\ref{con_2}),
which represents an example of a quantum-inspired scheme.

In such a classical LogQ formulation, the complexities associated with Pauli decomposition or quantum measurements are eliminated. LogQ can thus be efficiently implemented using classical computational frameworks.

\section{Conclusion}

In this paper we illustrated how quantum principles can inspire classical computation.
LogQ was originally developed in a quantum framework as a pragmatic approach to solve QUBO problems.
While recent work made quantum LogQ compliant with gradient-inspired schemes, the Pauli-decomposition and measurement overhead remained a practical challenge.

We presented a full reformulation of LogQ within a classical framework that bypasses this challenge while capitalizing on previously developed LogQ conditions. This classical variant functions as a continuous (non-linear) relaxation of variables acting on phases, a nod to its quantum origins. It can be used as a standalone heuristic or might be embedded within classical algorithms (like B\&B to compute incumbents, etc.). Further study would make sense to optimize the choice of the $R(\cdot)$ function (or even define other forms for $f(\cdot)$ that satisfy the first line of \eqref{eq:QUBO}) and to evaluate the scheme’s scalability on large-scale industrial instances where traditional QUBO solvers tend to struggle.

\section*{Acknowledgment}
The authors warmly thank Isabel Barros Garcia and Jean-Patrick Mascomère for insightful comments.
 The authors thank TotalEnergies for permission to publish this work.

\bibliography{conference_101719}

\end{document}